%
%

%
%


\magnification 1200
\input amstex
\documentstyle{amsppt}
\NoBlackBoxes
\NoRunningHeads


\define\qu{quasiunipotent}
\define\vre{\varepsilon}

\define\hs{homogeneous space}
\define\df{\overset\text{def}\to=}
\define\un#1#2{\underset\text{#1}\to#2}
\define\br{\Bbb R}
\define\bn{\Bbb N}
\define\bz{\Bbb Z}
\define\bq{\Bbb Q}

\define\di{Diophantine}
\define\da{Diophantine approximation}

\define\ve{\bold e}
\define\vx{\bold x}
\define\vy{\bold y}
\define\vz{\bold z}
\define\vv{\bold v}

\define\vw{\bold w}
\define\vp{\bold p}
\define\vq{\bold q}

\define\vc{\bold c}
\define\vf{\bold f}
\define\vt{\bold t}


\define\nz{\smallsetminus \{0\}}

\define\be{Besicovitch}

\define\cag{$(C,\alpha)$-good}

\define\thuss{$\Leftrightarrow$ }
\define\ca{\Cal A(d,k,C,\alpha,\rho)}
\define\sfb{(S,\varphi,B)}
\define\ph{\varphi}
\define\p{\Phi}
\define\vrn{\varnothing}
\define\ssm{\smallsetminus}
\define\vwa{very well approximable}

\topmatter
\title Flows on homogeneous spaces and Diophantine approximation on manifolds\endtitle  

\author { D.$\,$Y.~Kleinbock} \\ 
  { \rm 
   Institute for Advanced Study} \vskip 5pt 
 and G.$\,$A.~Margulis \\ 
  { \rm 
   Yale University} \\ \ 
\endauthor

    \address{Dmitry Y. Kleinbock, {\it Present address:}  Department of
Mathematics, Rutgers University, New Brunswick, NJ 08903
}
  \endaddress
\email kleinboc\@math.rutgers.edu
 \endemail

    \address{ G.$\,$A. Margulis: Department of Mathematics, Yale University, 
   New Haven, CT 06520
}
  \endaddress

\email  margulis\@math.yale.edu \endemail
  \thanks The work of the first author was supported in part by NSF Grant DMS-9304580 and that of the second author by NSF Grant DMS-9424613. \endthanks       

      
\abstract 
We present a new approach to  metric Diophantine approximation on manifolds based on the correspondence between approximation properties of numbers and orbit properties of certain flows on \hs s. This approach yields a new proof of a conjecture of Mahler, originally settled by V.$\,$G.~Sprind\v zuk in 1964. We also prove several related hypotheses  of Baker 
and Sprind\v zuk formulated in 1970s. The core of the proof is a theorem which generalizes and 
sharpens earlier results  on non-divergence of 
unipotent flows on the space of lattices.
 \endabstract


\dedicatory To appear in Annals of Mathematics
\enddedicatory

\endtopmatter
\document

\heading{1. Introduction}
\endheading 

We start by recalling several basic facts from the theory of simultaneous \da. For $\vx,\vy\in\br^n$ we let 
$$
\split
\vx\cdot\vy = \sum_{i = 1}^n x_iy_i&, \quad \|\vx\| = \max_{1\le i \le n}|x_i|,\\
\Pi(\vx) = \prod_{i = 1}^n |x_i|\quad&\text{ and }\quad\Pi_{\sssize +}(\vx) = \prod_{i = 1}^n |x_i|_{\sssize +},
\endsplit
$$
where $|x|_{\sssize +}$ stands for $\max(|x|, 1)$. 
 One says that a vector $\vy\in\br^n$ is {\it \vwa\/} (cf.~\cite{S2}), to be abbreviated as VWA,  if the following two equivalent conditions are satisfied:

(i) for some $\vre > 0$ there are infinitely many $\vq\in \bz^n$ such that 
$$
|\vq\cdot\vy + p|\cdot\|\vq\|^n \le \|\vq\|^{-n\vre}\tag 1.1
$$
for some $p\in\bz$;

(ii) for some $\vre > 0$ there are infinitely many $q\in \bz$ such that 
$$
\|q\vy + \vp\|^n \cdot |q| \le |q|^{-\vre}\tag 1.2
$$
for some $\vp\in\bz^{n}$. 

The equivalence of (1.1) and (1.2) follows from Khintchine's
Transference Principle, see e.g.~\cite{C, Chapter V}. It is well known
(see \cite{C, Chapter VII}) that almost every $\vy\in\br^n$ is not
VWA. A more difficult question arises if one considers almost all
points $\vy$ on a submanifold $M $ of $\br^n$ (in the sense of the
natural measure class on $M$). In 1932 K.~Mahler \cite{Mah} conjectured that almost all points on the curve 
 $$
M = \{(x,x^2,\dots,x^n)\mid x\in \br\}\tag 1.3
$$ 
are not very well approximable. V.~Sprind\v zuk's  proof of this conjecture (see \cite{Sp1, Sp2}) has eventually led to the development of a new branch in approximation theory, usually referred to as ``\da\ with dependent quantities'' or ``\da\ on manifolds''. 
According to Sprind\v zuk's terminology, a submanifold 
 $M  \subset \br^n$ is called {\it extremal\/} if almost all $\vy\in M $ are not VWA. Since 1964, extensive classes of extremal manifolds have been found by Sprind\v zuk and his followers, see \cite{{Sp3}, {Sp4}} for a review. However, the following conjecture made by Sprind\v zuk in 1980 has remained unproved:

\proclaim{Conjecture H${}_{\bold 1}$ \rm \cite{Sp4}} Let $f_1,\dots,f_n$ be real analytic functions in $\vx\in U$, $U$ a domain in $\br^d$,  which together with $1$ are linearly independent over $\br$. Then the manifold
$
M  = \{\vf(\vx)\mid \vx\in U\}
$
 is extremal.
\endproclaim

The case $n = 2$ was settled by  Schmidt in 1964 \cite{S1}, and very recently    V.~Beresnevich and V.~Bernik \cite{BeBe} proved the above result for 
$n = 3$. 

Note that the validity of the above conjecture for polynomials
$f_1,\dots,f_n \in \bq[x]$  
with deg$(f_i) \le n$  follows from Sprind\v zuk's original proof of
Mahler's Conjecture. However, there exists a strengthening of
Conjecture H${}_1$ which has not been proved even for the curve
(1.3). Namely, one says that $\vy\in\br^n$ is {\it very well
multiplicatively approximable\/}  (to be abbreviated as VWMA) if the
following two equivalent conditions are satisfied: 

(i) for some $\vre > 0$ there are infinitely many $\vq\in \bz^n$ such that 
$$
|\vq\cdot\vy + p|\cdot\Pi_{\sssize +}(\vq) \le \Pi_{\sssize +}(\vq)^{-\vre}\tag 1.1M
$$
for some $p\in\bz$;

(ii) for some $\vre > 0$ there are infinitely many $q\in \bz$ such that 
$$
\Pi(q\vy + \vp)\cdot |q| \le |q|^{-\vre}\tag 1.2M
$$
for some $\vp\in\bz^{n}$. 

In other words, instead of taking the $n$th power of the maximum of the ``differences'' $|qy_i + p_i|$ in (1.2) one simply multiplies them, and instead of taking the $n$th power of the maximum of the numbers $|q_i|$ in (1.1) one multiplies those of them which are different from zero. The equivalence of (1.1M) and (1.2M) can be shown by modification of the argument needed to prove the standard form of  Khintchine's Transference Principle, see a remark in \cite{Sp3, p.~69}.  A manifold 
 $M  \subset \br^n$ is said to be {\it strongly extremal\/} if almost all $\vy\in M$ are not VWMA. Clearly (1.$i$) implies (1.$i$M), $i = 1,2$; therefore any strongly extremal manifold is extremal. In his book \cite{B}, A.~Baker raised the question of confirming the strong extremality of the curve (1.3). Later, in his review \cite{Sp4}, Sprind\v zuk also stated 

\proclaim{Conjecture H${}_{\bold 2}$}  Any manifold $M \subset \br^n$ satisfying the assumptions of Conjecture {\rm H}${}_1$ is strongly extremal.
\endproclaim

We remark that the validity of the above conjecture for polynomials with $n \le 4$ is the recent result of V.~Bernik and V.~Borbat  \cite{BeBo}. See also \cite{VC} and references therein for other special cases.

In this paper we present a proof of the above conjectures as well as some of their generalizations. More precisely, consider   a $d$-dimensional submanifold $
M  = \{\vf(\vx)\mid \vx\in U\}
$ of $\br^n$, where $U$ is an open subset of $\br^d$ and $\vf = (f_1,\dots,f_n)$ is a $C^m$ imbedding of $U$ into $\br^n$. For $l\le m$, say that $\vy = \vf(\vx)$ is an {\it $l$-nondegenerate point\/} of $M$ if   the space $\br^n$ is spanned by partial derivatives of $\vf$ at $\vx$ of order up to $l$. We will say that $\vy$ is 
 {\it nondegenerate\/} if it is $l$-nondegenerate for some $l$. One can view this condition as an infinitesimal version of not lying in any proper affine hyperplane, i.e.~of the linear independence of $1,f_1,\dots,f_n$  over $\br$. Indeed, if the functions $f_i$ are analytic, it is easy to see that the linear independence of $1,f_1,\dots,f_n$  over $\br$ in a domain $U$ is equivalent to all points of $M = \vf(U)$ being {nondegenerate}. Thus {Conjecture H${}_{2}$} would follow from

\proclaim{Theorem A} Let $f_1,\dots,f_n\in C^m(U)$, $U$ an open subset of $\br^d$,  be such that a.e.~point of $M = \{\vf(\vx)\mid \vx\in U\}$ is {nondegenerate}.
Then $M$
 is strongly extremal.
\endproclaim

Note that if $n = 2$ and $d = 1$, the above assumption says that the curve $\big\{\big(f_1(x),f_2(x)\big)\mid x\in U\big\}$, $U\subset \br$, has nonzero curvature almost everywhere. This was exactly the assumption used by Schmidt \cite{S1}, and the strong extremality of such a curve is the subject of Theorem 14 from \cite{Sp3, Chapter 2}.

 Our proof of Theorem A is based on the correspondence (cf.~\cite{D1, K1, K2}) between approximation properties of numbers $\vy\in \br^n$ and behavior of certain orbits in the space of unimodular lattices in $\br^{n+1}$. More precisely, one acts on the lattice
$$
\Lambda_\vy = \left(\matrix
1 & \vy^{\sssize T}  \\
0 & I_n
\endmatrix \right)\bz^{n+1}\tag 1.4
$$
by semisimple elements of the form
$$
g_\vt = \text{\rm diag}(e^{t},e^{-t_1},\dots,e^{-t_n}), \ \vt = (t_1,\dots,t_n),\ t_i \ge 0,\ t = \sum_{i = 1}^nt_i\tag 1.5
$$
(the latter notation will be used throughout the paper, so that whenever $t$ and $\vt$ appear in the same formula, $t$ will stand for $\sum_{i = 1}^nt_i$).

Define  a function $\delta$ on the space of lattices by  
$$
\delta(\Lambda) \df \inf_{\vv\in\Lambda\nz}\|\vv\|\tag 1.6
$$
(note that the ratio of $1 + \log\big(1/{\delta(\Lambda)}\big)$ and $1 + \text{dist}(\Lambda, \bz^{n+1})$ is bounded between two positive constants for any right invariant Riemannian metric ``dist'' on the space of lattices). We    prove in \S 2 that  
for any  very well multiplicatively approximable $\vy\in\br^n$ 
  there exists $\gamma > 0$ and infinitely many 
$\vt\in\bz_{\sssize +}^n$ such that
$$
\delta(g_\vt\Lambda_{\vy}) \le e^{-\gamma t}\,;\tag 1.7
$$
in other words, if dist$(g_\vt\Lambda_\vy, \bz^{n+1})$ grows sublinearly as a function of $\vt$, then $\vy$ is not VWMA. 
Thus to prove the strong extremality of $
M$ as in Theorem A it is enough to show that  for any nondegenerate point $\vy_0 = \vf(\vx_0)$ there is a neighborhood $B$ of $\vx_0$ in $U$ such that for almost all $\vy\in \vf(B)$ and any  $\gamma > 0$, there are at most finitely many $\vt\in\bz_{\sssize +}^n$ such that (1.7) holds.

In view of
Borel-Cantelli, the latter can be done by estimating the measure of the sets 
$$
E_\vt\df\{\vx\in B\mid \delta(g_\vt\Lambda_{\vf(\vx)}) \le e^{-\gamma t}\}\tag 1.8
$$
 for any fixed $\vt\in \bz_{\sssize +}^n$, so that 
$$
\sum_{\vt\in\bz_{\sssize +}^n} |E_\vt|  < \infty\tag 1.9
$$
(here and hereafter  $|\cdot|$ stands for the Lebesgue measure).
Such estimates are obtained in \S 5 by modifying proofs of earlier results  on non-divergence of unipotent flows in the space of lattices. According to the theorem of Dani \cite{D2} strengthening an earlier result of Margulis \cite{Mar},
for any $c > 0$ and any lattice $\Lambda$ in $\br^k$ there exists $\vre > 0$ such that  for any unipotent subgroup $\{u_x|x\in\br\}$ of $SL_k(\br)$ one has
$$
\big|\{x\in[0,T]\mid \delta(u_x\Lambda) < \vre\}\big| \le cT \,.\tag 1.10
$$
Similar estimates are known for any polynomial map from $\br^d$ to $GL_k(\br)$ instead of $x\to u_x$, cf.~\cite{Sh} (or \cite{EMS} for a bigger class of maps). In \S 3 (see Proposition 3.4 for a precise statement) we show that if $\vy_0 = \vf(\vx_0)$ is nondegenerate, $\vx_0$  has a neighborhood  on which linear combinations of  $1,f_1,\dots,f_n$ behave like polynomials of uniformly bounded degree. Then in \S\S\ 4 and 5 we modify the argument of Margulis and Dani in order  to get a quantitative relation between $c$ and $\vre$ in the analogue of (1.10) (see Proposition 2.3) which will guarantee convergence in (1.9). 
The last section of the paper deals with  several possible extensions of the main result, as well as some open questions. 

\heading{2. Reduction to a statement about lattices}
\endheading 

Given a vector $\vy\in\br^n$, consider a lattice $\Lambda_\vy$ in $\br^{n+1}$ defined as in (1.4). Note that elements of $\Lambda_\vy$ are of the form 
$$
(\vq\cdot\vy + p,\ q_1,\dots,q_n)\,,\tag 2.1
$$
where $p\in\bz$ and 
$\vq = (q_1,\dots,q_n)\in \bz^n$. We are going to consider the action of the semigroup $\{g_\vt\}$ as in (1.5) on the space  of unimodular lattices in $\br^{n+1}$, i.e.~the \hs\ $SL_{n+1}(\br)/SL_{n+1}(\bz)$. The function $\delta$ defined in (1.6) will be used to ``describe the structure of the space  of  lattices  at $\infty$'': by Mahler's Compactness Criterion (see \cite{R, Corollary 10.9}) a sequence  $\Lambda_k$ tends to infinity in $SL_{n+1}(\br)/SL_{n+1}(\bz)$ iff  $
\delta(\Lambda_k)\to 0$ as $k\to\infty$.

The next lemma helps one to reduce problems involving approximation properties of $\vy$ to studying the rate, calculated in terms of the function $\delta$, with which orbit points $\{g_\vt\Lambda_\vy\}$ may go to infinity.

\proclaim{Lemma 2.1} Let $\vre > 0$, $\vy\in\br^{n}$ and $(p,\vq)\in\bz^{n+1}$ be such that {\rm (1.1M)} holds. Put 
$$
r = \Pi_{\sssize +}(\vq)^{-\frac{\vre}{n+1}}\,,\tag 2.2a
$$
 and for $i = 1,\dots,n$ define  $t_i> 0$ by 
$$
|q_i|_{\sssize +} = r e^{t_i}\,.\tag 2.2b
$$
Then 
$
\delta(g_\vt\Lambda_\vy) \le r
$.
\endproclaim

\demo{Proof} In view of (2.1) and (1.5), the inequality to be proved
would follow from the inequalities
$$
  e^{t}|\vq\cdot\vy + p| \le r\tag 2.3a
$$
and
$$
e^{-t_i}|q_i|\le r,\ i = 1,\dots,n\,.\tag 2.3b
$$
One immediately deduces (2.3b) from (2.2b) and $|q_i|\le
|q_i|_{\sssize +}$ Taking a product of the equalities (2.2b), one gets  
$$
\Pi_{\sssize +}(\vq) = r^ne^{t}\,.\tag 2.4
$$
Thus (1.1M) can be written in the form
$$ 
 e^{t}|\vq\cdot\vy + p| \le r^{-n}\Pi_{\sssize +}(\vq)^{-\vre} \underset{\text{by (2.2a)}}\to =  r^{-n}\cdot r^{n+1}\,,
$$
which proves (2.3a).
 \qed
\enddemo

An elementary computation shows  that (2.2a) and (2.4) can be combined to yield
$r = e^{-\gamma t}$, where  
$$
\gamma = \frac\vre{n + 1 + n\vre}\,.\tag 2.5 
$$
Using this remark, we prove

\proclaim{Corollary 2.2} Assume that $\vy\in\br^n$ is VWMA. Then  
  for some $\gamma > 0$ there are infinitely many 
$\vt\in\bz_{\sssize +}^n$ such that {\rm (1.7)} holds.
\endproclaim

\demo{Proof} By definition, for some $\vre> 0$ there are infinitely many solutions $\vq\in\bz^n$ of (1.1M). Therefore, by the above lemma and with $\gamma$ as in (2.5), there exists a sequence $\vt^{(k)}\in\br_{\sssize +}^n$ with $t^{(k)} = \sum_{i = 1}^nt_i^{(k)}\to + \infty$ such that $
\delta(g_{\vt^{(k)}}\Lambda_\vy) \le e^{-\gamma t^{(k)}}
$. \linebreak Denote by $[\vt^{(k)}]$ the vector consisting of  integer parts of $t_i^{(k)}$, then clearly the ratio of $\delta(g_{\vt^{(k)}}\Lambda_\vy)$ and $\delta(g_{[\vt^{(k)}]}\Lambda_\vy)$ is bounded from above by $$
\|g_{\vt^{(k)}}g^{-1}_{[\vt^{(k)}]}\| = \|g_{\vt^{(k)} - [\vt^{(k)}]}\| \le e^n\,.
$$
 Thus one gets  infinitely many solutions $\vt\in\bz_{\sssize +}^n$ of  $
\delta(g_\vt\Lambda_\vy) \le e^ne^{-\gamma t}$, and a slight variation of $\gamma$ yields  infinitely many solutions of (1.7). \qed
\enddemo

Now we are ready to start the

\demo{Proof of Theorem A} We are given an open subset $U$ of $\br^d$ and a smooth manifold $M$ parametrized by $\vf:U\to\br^n$. It is  enough to prove that almost every point of $M$ has a neighborhood $W\subset M$ such that the set $\{\vy\in W\mid \vy \text{ is not VWMA}\}$ has full measure. By virtue of Corollary 2.2, it suffices to show that
for a.e.~$\vy_0\in M$ there exists a neighborhood $W\subset M$ of $\vy_0$ such that for any $\gamma > 0$ and a.e.~$\vy\in W$, the inequality (1.7) holds for at most finitely many $\vt\in\bz_{\sssize +}^n$.

We now state a proposition, to be proved in \S 5, which will easily imply the above statement.

\proclaim{Proposition 2.3} Let $\vf$ be a $C^l$ map from an open subset $U$ of $\br^d$ to $\br^n$, and let $\vx_0\in U$  be such that $\br^n$ is spanned by partial derivatives of $\vf$ at $\vx_0$ of order up to $l$. Then there exists a ball $B\subset U$ centered in $\vx_0$ and positive constants $D$ and $\rho$ such that  for any $t_1,\dots,t_n \ge 0$ and  $0 < \vre \le \rho$ one has
$$
\big|\{\vx\in B\mid \delta\big(g_\vt \Lambda_{\vf(\vx)}\big) < \vre\}\big| \le D\left(\frac\vre\rho\right)^{1/dl}|B|\,.\tag 2.6
$$
\endproclaim

Indeed, by the assumptions on $M$, almost every $\vy_0\in M$ is nondegenerate, hence has a neighborhood $W = \vf(B)$ with $B$ as in the above proposition. For $\gamma > 0$ and  $\vt\in\bz_{\sssize +}^n$ consider the set $E_\vt$ defined in (1.8). Then there is a constant $\tilde D$ independent on $\vt$ such that for large enough $\|\vt\|$ one has  
$
|E_\vt| \le \tilde D e^{-\gamma t/dl}
$. Therefore the series $
\sum_{\vt\in\bz_{\sssize +}^n} |E_\vt|$ converges, hence, by the Borel-Cantelli Lemma, for  a.e.~$\vx\in B$ the inequality (1.7), with $\vy = \vf(\vx)$, holds  for at most finitely many $\vt\in\bz_{\sssize +}^n$.
Since the measure on $M$ is obtained from the Lebesgue measure on $U$ by a smooth map, this finishes the proof of Theorem A modulo Proposition 2.3. \qed
\enddemo

\heading{3. Good functions and nondegenerate points}
\endheading

Let  $V$ be a subset of $\br^d$ and $f$ a continuous function on $V$. In what follows, we will let $\|f\|_{\sssize B} \df \sup_{\vx\in B}|f(\vx)|$ for a subset $B$ of $V$. For  positive numbers $C$ and $\alpha$, say that 
$f$ is {\it \cag\ on $V$\/} if for any open ball $B\subset V$
one has \footnote{Here we adopt the convention $\frac10 = \infty$, so that (3.1) holds if $f|_{\sssize B}\equiv 0$.}
$$
\forall\,\vre > 0\quad |\{\vx\in B\bigm| |f(\vx)| < \vre\}| \le C\cdot\left(\frac\vre{\|f\|_{\sssize B}}\right)^\alpha \cdot |B|\,.
\tag 3.1
$$
 
Cf.~\cite{EMS, Corollary 2.10} for a similar condition.   
The properties listed below follow immediately from the definition.

\proclaim{Lemma 3.1} Let $V\subset \br^d$ and $C,\alpha > 0$ be given.

{\rm (a)} $f$ is \cag\ on $V$ \thuss so is $|f|$;

{\rm (b)} $f$ is \cag\ on $V$ $\Rightarrow$  so is $\lambda f$ $\forall\,\lambda\in \br$;

{\rm (c)}  $f_i$, $i\in I$, are \cag\ on $V$ $\Rightarrow$ so is $\sup_{i\in I}|f_i|$;
\endproclaim


We now consider what can be called a model example of good functions.

\proclaim{Proposition 3.2 \rm (cf.~\cite{DM, Lemma 4.1})} For any $k\in\bn$, any polynomial $f\in\br[x]$ of degree not greater than $k$ is $\big(2k(k+1)^{1/k}, 1/k\big)$-good on $\br$.
\endproclaim

\demo{Proof} Fix an open interval $B\subset \br$, a  polynomial $f\in\br[x]$ of degree not exceeding $k$ and an $\vre > 0$. We need to show that
$$
\left|\{x\in B\bigm| |f(x)| < \vre\}\right| \le 2k(k+1)^{1/k}\left(\frac\vre{\|f\|_{\sssize B}}\right)^{1/k}|B|\,.\tag 3.2
$$

Denote $|B|$ by $b$ and $\dfrac1b\left|\{x\in B\bigm| |f(x)| < \vre\}\right|$ by $\sigma$. Then there exist $x_1,\dots,x_{k+1} \in B$ such that $|f(x_i)| \le \vre$, $1\le i \le k+1$, and $|x_i - x_j| \ge \sigma b/2k$, $1\le i < j \le k+1$. Using Lagrange's interpolation formula, one can write down the exact expression for $f$:
$$
f(x) = \sum_{i = 1}^{k+1}f(x_i)\frac{\prod_{j = 1,\,j\ne i}^{k+1}(x - x_j)}{\prod_{j = 1,\,j\ne i}^{k+1}(x_i - x_j)}\,,\tag 3.3a
$$
and conclude that
$$
\|f\|_{\sssize B} \le (k + 1)\vre \frac{b^k}{(\sigma b/2k)^k}\,,\tag 3.3b
$$
which immediately implies (3.2). \qed \enddemo

The goal of this section is to prove that if the point
$\vy_0 = \vf(\vx_0)$ of the manifold $M = \vf(U)\subset \br^n$, $U\subset \br^d$, is nondegenerate, then there exist a neighborhood $V$ of $\vx_0$ and positive constants $C$ and $\alpha$ such that all linear combinations of $1,f_1,\dots, f_n$ are \cag\ on $V$. Our  argument is based on the following standard but quite technical lemma.
 For $i = 1,\dots,d$, we let $\partial_i$ be the operator of partial differentiation with respect to $x_i$, and for a {\it multiindex\/} $\beta = (i_1,\dots,i_d)$, $i_j\in \bz_{\sssize +}$, we let $|\beta| = \sum_{j = 1}^d i_j$ and $\partial_\beta = \partial_1^{i_1}\circ\dots\circ \partial_d^{i_d}$.

\proclaim{Lemma 3.3} Let $V$ be an open subset of $\br^d$,  and let $f\in C^k(V)$ be such that for some constants $A_1, A_2 > 0$ one has 
$$
\|\partial_\beta f\|_{\sssize V} \le A_1 \quad \forall\,\beta\text{ with }|\beta|\le k\,,\tag 3.4$\le$
$$
and 
$$
|\partial_i^{k}f(x)|\ge A_2\quad \forall\,x\in V,\ i = 1,\dots,d\,.\tag 3.4$\ge$
$$
Then for any $d$-dimensional cube $B = B_1\times\dots\times B_d\subset V$, where $B_i$ are open intervals of the same length, and for any $\vre > 0$ one has
$$
\left|\{\vx\in B\bigm| |f(\vx)| < \vre\}\right| \le d C_{\sssize k,A_1,A_2} \left(\frac\vre{\|f\|_{\sssize B}}\right)^{1/dk}|B|\,,\tag 3.5
$$
with 
$$
C_{\sssize k,A_1,A_2} = k(k+1)\left(\frac{A_1}{A_2}(k+1)(2k^k+1)\right)^{1/k}\,.
$$
\endproclaim

\demo{Proof} For $B$, $f$ as above and $\vre > 0$, we will denote  by  $b$ the sidelength of $B$ and by $s$ the supremum of $|f|$ over $B$. First consider the case $d = 1$. Take $\vre > 0$; since, by (3.4$\ge$),  the $k$th derivative of $f$ does not vanish on $B$, the set $\{x\in B\bigm| |f(x)| < \vre\}$ consists of at most $k + 1$ intervals. Let $I$ be the maximal of those, then  
$$
\left|\{x\in B\bigm| |f(x)| < \vre\}\right| \le (k + 1)|I|\,,\tag 3.6
$$
 so it suffices to estimate $|I|$ from above.  Denote by $r$ the length of $I$ and also let 
$$
\tau  = k!{s}/{b^k}\,.\tag 3.7
$$
The argument will be based on the following two inequalities: 

\proclaim{Sublemma 3.3.1} {\rm (a)} $ r \le \ph_1(\tau )  \left(\dfrac\vre{s}\right)^{1/k}b$, where $\ph_1(\tau ) = k\left(\dfrac{k+1}{A_2}\right)^{1/k}\tau ^{1/k}$;

{\rm (b)}  $ r \le \ph_2(\tau )\left(\dfrac\vre{s}\right)^{1/k}b$, where $\ph_2(\tau ) =  k^2\left(\dfrac{\frac{2 A_1}{A_2}(k+1)}{1 - {A_1}/{\tau }}\right)^{1/k}$. \endproclaim

\demo{Proof} Divide $I$ into $k$ equal parts by points $x_1,\dots,x_{k+1}$, and let $P$ be the Lagrange polynomial of degree $k$ formed by using values of $f$ at these points, i.e.~given by the expression in the right hand side of (3.3a). Then there exists $x\in I$ such that $P^{(k)}(x) = f^{(k)}(x)$, hence, by (3.4$\ge$), $|P^{(k)}(x)| \ge A_2$. On the other hand, one can differentiate the right hand side of (3.3a) $k$ times to get  
$|P^{(k)}(x)| \le (k + 1) \dfrac{\vre k!}{(r/k)^k}$. Combining the last two inequalities, one obtains
$$
 r \le 
k\left(\dfrac{(k+1)!}{A_2}\vre\right)^{1/k} = k\left(\dfrac{(k+1)!}{A_2}\dfrac{s}{b^k}\right)^{1/k}\left(\dfrac\vre{s}\right)^{1/k}b\,,\tag 3.8
$$
which reduces to the inequality (a) by the substitution (3.7).

Next, let $Q$ be the Taylor polynomial of $f$ at $x_1$ of degree $k-1$. By Taylor's formula, 
$$
\|f - Q\|_{\sssize I} \le \|f^{(k)}\|_{\sssize I} \frac{r^k}{k!} \un{by (3.4$\le$)}\le  \frac{A_1r^k}{k!} \un{by (3.8)}\le \frac{A_1}{A_2}k^k(k+1)\vre\,,
$$ 
therefore
$
\|Q\|_{\sssize I} \le 2\dfrac{A_1}{A_2}k^k(k+1)\vre
$.

We now apply Lagrange's formula to reconstruct $Q$ on $B$ by its values at $x_1,\dots,x_{k}$. Similarly to (3.3b), we get
$$
\|Q\|_{\sssize B} \le k\cdot 2\frac{A_1}{A_2}k^k(k+1)\vre \frac{b^{k-1}}{(r/k)^{k-1}}\le 2\frac{A_1}{A_2}k^{2k}(k+1)\vre \frac{b^k}{r^k}\,.\tag 3.9
$$
Finally, the difference between $f$ and $Q$ on $B$ is, again by (3.4$\le$), bounded from above by ${A_1b^k}/{k!}$, so from (3.9) one deduces that 
$$
s \le \left(2\frac{A_1}{A_2}k^{2k}(k+1) \frac{\vre}{r^k} +\frac{A_1}{k!}\right)b^k\,,
$$
which is equivalent to the inequality (b) modulo (3.7).
\qed\enddemo

To use the above estimates, note that $\ph_1$ (resp.~$\ph_2$) is an increasing (resp.~decreasing) function, therefore one has $ r \le \ph_1(\tau_0)  \left(\dfrac\vre{s}\right)^{1/k}b$, where $\tau_0$ is the root of $ \ph_1(\tau) = \ph_2(\tau)$. An elementary computation yields $\tau_0 = A_1(2k^k + 1)$, and the validity of (3.5) for $d = 1$  follows immediately from (3.6). 

Next we argue by induction on $d$, assuming that (3.5) is established for all the lower values of  $d$. We represent $\vx\in\br^d$ as $(x_1,\vx')$ and let $B'$ be the product of the intervals $B_2,\dots,B_d$. Denote by $\sigma$ the value $\dfrac1{ b^d}\left|\{\vx\in B\bigm| |f(\vx)| < \vre\}\right|$. Choose a positive number $\lambda < \sigma$, and let
$$
B_1(\lambda)\df \left\{x_1\in B_1\left| \frac1{b^{d-1}}\right.\left|\{\vx'\in B'\bigm |f(x_1,\vx')| < \vre\}\right| \ge \lambda\right\}\,.
$$

The induction step will be based on the following 

\proclaim{Sublemma 3.3.2} $\dfrac1{b}|B_1(\lambda)| \ge \sigma - \lambda$. \endproclaim

\demo{Proof} Write
$$
\{\vx\in B\bigm| |f(\vx)| < \vre\}\subset B_1(\lambda)\times B' \ \cup \ \left\{\vx\in \big(B_1\ssm B_1(\lambda)\big)\times B'\bigm| |f(\vx)| < \vre\right\}\,,
$$
and
put $\mu = \dfrac1{b}|B_1(\lambda)|$. Then, by Fubini, $b^d\sigma = \left|\{\vx\in B\bigm| |f(\vx)| < \vre\}\right|$ is not greater than $b\mu\cdot b^{d-1} + b(1 - \mu)\cdot b^{d-1}\lambda \le b^d(\lambda + \mu)$, hence the claim. \qed\enddemo

For $x_1\in B_1$, denote by $f_{x_1}$ the function on $B'$ given by $f_{x_1}(\vx') \df f(x_1,\vx')$. Clearly such functions satisfy assumptions (3.4) on derivatives with $d -1$ in place of $d$, therefore, by the induction assumption, for any $x_1\in B_1(\lambda)$ one has
$$
\lambda \le \frac1{b^{d-1}}\left|\left\{\vx'\in B'\bigm |f_{x_1}(\vx')| < \vre\right\}\right| \le (d-1) C_{\sssize k,A_1,A_2} \left(\frac\vre{\|f_{x_1}\|_{\sssize B'}}\right)^{1/(d-1)k}\,,
$$
or
$$
\|f_{x_1}\|_{\sssize B'} \le \left(\frac{(d-1)C_{\sssize k,A_1,A_2}}\lambda\right)^{(d-1)k} \vre\,.
$$
On the other hand, for any $\vx'\in B'$ the functions $x_1\to f(x_1,\vx')$  also satisfy the assumptions on derivatives (with $1$ in place of $d$), hence are $\left(C_{\sssize k,A_1,A_2}, 1/k\right)$-good on $B_1$. By Lemma 3.1(c), the function $x_1\to \|f_{x_1}\|_{\sssize B'}$ is also  $\left(C_{\sssize k,A_1,A_2}, 1/k\right)$-good on $B_1$, thus, by the above sublemma,
$$
\split
\sigma - \lambda &\le \frac1{b}\left|\left\{x_1\in B_1\left|\|f_{x_1}\|_{\sssize B'} \le \big({(d-1)C_{\sssize k,A_1,A_2}}/\lambda\big)\right. ^{(d-1)k} \vre\right\}\right|\\
& \le C_{\sssize k,A_1,A_2}\left(\frac{\big({(d-1)C_{\sssize k,A_1,A_2}}/\lambda\big)^{(d-1)k} \vre}s\right)^{1/k}\,,
\endsplit
$$
or equivalently,
$$
\lambda^{d-1}(\sigma - \lambda) \le C_{\sssize k,A_1,A_2}^d  (d-1)^{d-1}\left(\frac{\vre}s\right)^{1/k}\,.\tag 3.10
$$
The function in the left hand side of (3.10) attains its maximum, $\sigma^d{(d-1)^{d-1}}/{d^{d}}$, when $\lambda = \sigma({d-1})/{d}$; substituting it into (3.10), one easily obtains (3.5).
\qed\enddemo

We now come to the main result of the section.

\proclaim{Proposition 3.4} Let $\vf = (f_1,\dots,f_n)$ be a $C^l$ map from an open subset $U$ of $\br^d$ to $\br^n$, and let $\vx_0\in U$ be such that $\br^n$ is spanned by partial derivatives of $\vf$ at $\vx_0$ of order up to $l$. Then there exists a neighborhood $V\subset U$ of $\vx_0$ and positive $C$  such that any linear combination of $ 1,f_1,\dots,f_n$ is $(C,1/dl)$-good on $V$. \endproclaim

\demo{Proof} Take $f = c_0 + \sum_{i = 1}^nc_if_i$; in view of Lemma 3.1(b), one can without loss of generality assume that the norm of $\vc = (c_0,\dots,c_n)$ is equal to $1$.  
From the nondegeneracy assumption it follows that there exists a constant $C_1 > 0$ such that for any $\vc$ with $\|\vc\| = 1$ one can find a multiindex $\beta$ with $|\beta| = k \le l$ and  
$$
\big|\sum_{i = 1}^nc_i\partial_\beta f_i(\vx_0)\big| = |\partial_\beta f(\vx_0)|\ge C_1\,.
$$
 By an appropriate rotation  of the coordinate system around $\vx_0$ one can guarantee that 
$|\partial_i^k f(\vx_0)|\ge C_2$ for  all $i = 1,\dots,d$ and some positive $C_2$ independent of $\vc$. Then one uses the continuity of the derivatives of $f_1,\dots,f_n$ to choose a neighborhood $V'\subset U$ of $\vx_0$ and positive $A_1$, $A_2$  (again independently of $\vc$) such that the inequalities (3.4) hold. Now let $V$ be a smaller neighborhood of $\vx_0$ such that whenever a ball $B$ lies in $V$, any cube $\hat B$ circumscribed around $B$ is contained in $V'$. Then  for any $\vre > 0$ one has
$$
\split
\left|\{\vx\in B\bigm| |f(\vx)| < \vre\}\right| &\le \left|\{\vx\in \hat B\bigm| |f(\vx)| < \vre\}\right|\\ \text{(by Lemma 3.3)}\quad&\le  d C_{\sssize k,A_1,A_2} \left(\frac\vre{\|f\|_{\sssize \hat B}}\right)^{1/dk}|\hat B| \le  \frac{2^d}{v_d}d C_{\sssize k,A_1,A_2} \left(\frac\vre{\|f\|_{\sssize B}}\right)^{1/dk}|B|
\endsplit
$$
(here $v_d$ stands for the volume of the unit ball in $\br^d$), 
which implies that $f$ is $\left(\frac{2^d}{v_d}d C_{\sssize l,A_1,A_2}, 1/{dl}\right)$-good on $V$. 
\qed\enddemo

\heading{4. Maps of posets into spaces of good functions}
\endheading 

In this section we will work with  mappings of partially ordered sets ({\it
posets\/}) into spaces of functions on balls in $\br^d$. Given  a
mapping  from a poset  to the space of functions on a ball
$B$, we will {\it mark\/} certain points of $B$ (see the definition
below), and prove an upper
estimate (Theorem 4.1) for  the measure of the set of
unmarked points. Then in \S 5 we will use this estimate to  generalize
and strengthen results on  
non-divergence of unipotent flows on spaces of lattices obtained in
\cite{Mar, D2}. 

 In what follows, $B(\vx,r)$, where $\vx\in\br^d$ and $r > 0$, will stand for the open ball of radius $r$ centered in $\vx$. 
 For a 
poset $S$, we will denote by $l(S)$ the {\it length\/} of $S$ (i.e.~the number of elements in a maximal linearly ordered subset of $S$). If $T$ is a subset of $S$, we let $S(T)$ be the poset of elements of $S\smallsetminus T$ comparable with any element of $T$.
Note that one always has
$$
l\big(S(T)\big) \le l(S) - l(T)\,.\tag 4.1
$$

For $d\in \bn$, $k\in \bz_{\sssize +}$ and $C,\alpha,\rho   > 0$, define $\ca$ to be the set of triples $\sfb$ where $S$ is a poset,  $B = B(\vx_0,r_0)$, where $\vx_0\in\br^d$ and $r_0 > 0$,  and $\ph$ is a mapping from $S$ to the space of continuous functions on $\tilde B \df B\big(\vx_0,3^kr_0\big)$ (this mapping will be denoted by $s\to \psi_s$) such that the following holds:

\roster
\item"(A0)" $l(S) \le k$;

\item"(A1)" $\forall\,s\in S\,,\quad \psi_s$ is \cag\ on $\tilde B$;

\item"(A2)" $\forall\,s\in S\,,\quad\|\psi_s\|_{\sssize B} \ge \rho $;

\item"(A3)"  $\forall\,\vx\in \tilde B,\quad\#\{s\in S\bigm| |\psi_s(\vx)| < \rho\} < \infty$.
\endroster

Then, given $\sfb\in\ca$ and $\vre > 0$, say that a point $\vz\in B$ is {\it $(\vre,S,\ph)$-marked\/} if there exists a linearly ordered subset $\Sigma_\vz$ of $S$ such that

\roster
\item"(M1)" $\vre\le |\psi_s(\vz)| \le \rho \quad \forall\,s\in\Sigma_\vz$;

\item"(M2)" $|\psi_s(\vz)| \ge \rho \quad\forall\,s\in S(\Sigma_\vz)$.
\endroster

We will denote by $\p(\vre,S,\ph,B)$ the set of all the $(\vre,S,\ph)$-marked points $\vz\in B$. 
We will also need to use

\proclaim{Besicovitch's Covering Theorem \rm (see \cite{Mat, Theorem 2.7})}   There is an integer $N_d$ depending only on $d$ with the following property: let  $A$ be a bounded subset of $\br^d$ and let $\Cal B$ be a family  of nonempty open balls in $\br^d$ such that each $\vx\in A \text{ is the center of some ball of }\Cal B$; 
then there exists a finite or countable subfamily $\{B_i\}$ of $\Cal B$ with $
1_A \le \sum_i 1_{B_i} \le N_{d}
$ (i.e.~$A\subset \bigcup_i B_i$ and the multiplicity of that
subcovering is at most $N_{d}$). 
\endproclaim

The goal of the section is to prove the following

 \proclaim{Theorem 4.1 \rm (cf.~\cite{Mar, Main Lemma} or \cite{D2, Proposition 2.7})} Let $d\in \bn$, $k\in \bz_{\sssize +}$ and $C,\alpha,\rho  > 0$ be given. Then  for all $\sfb\in\ca$ and $ \vre >0$ one has
$$
\left|\big(B\smallsetminus \p(\vre,S,\ph,B)\big)\right| \le kC\big(3^dN_d\big)^k \left(\frac\vre \rho \right)^\alpha  |B|\,.
$$
\endproclaim
 
The proof will be built up from Lemmas 4.2--4.6 below. In these lemmas,  $d$, $k$, $C,\alpha$ and $\rho$ will be as in the  above theorem, and   $\sfb\in\ca$ will be fixed. We also
 define 
$$
H(\vx)\df\{s\in S\bigm| |\psi_s(\vx)| < \rho \}
$$
for any  $\vx\in B$ (this is a finite subset of $S$ in view of (A3)), and let 
$$
E\df \{\vx\in B\mid H(\vx)\ne\vrn\} = \{\vx\in B\mid \exists\,s\in S\text{ with }|\psi_s(\vx)| < \rho \}\,.
$$

\proclaim{Lemma 4.2} Any point of $B$ which does not belong to $E$ is $(\vre,S,\ph)$-marked for any  positive $\vre$. In other words, 
$
B\ssm E \subset \p(\vre,S,\ph,B)
$, and therefore
$$
B\ssm \p(\vre,S,\ph,B)\subset E\ssm \p(\vre,S,\ph,B)\,.
$$
\endproclaim

\demo{Proof} Take $\vz\in B\ssm E$; by the definition of $E$, one has $|\psi_s(\vz)| \ge \rho $ for all $s\in S$. Then one can take $\Sigma_\vz$ to be the empty set and check that (M1) and (M2) are satisfied. \qed\enddemo

The next four lemmas deal with some properties of the set $E$. Note that if $l(S)= 0$ (which means $S = \vrn$), one has $H(\vx)= \vrn$ for all $\vx\in B$, therefore $E = \vrn$. Thus in the argument below  we will tacitly assume $k$ to be not less than $1$. 

Take $\vx\in E$ and $s\in H(\vx)$, and define
$$
r_{s,\vx}\df\sup\{0< r \le 2r_0\bigm| \|\psi_s\|_{\sssize B(\vx,r)} \le \rho \}\,.
$$
 It follows from the continuity of functions $\psi_s$ that 
$r_{s,\vx} > 0$.  
We also  let 
$B_{s,\vx}\df B(\vx,r_{s,\vx})$.

\proclaim{Lemma 4.3}   For any $\vx\in E$ and $s\in H(\vx)$, one has 
$\|\psi_s\|_{\sssize B_{s,\vx}} \ge \rho$.
\endproclaim

\demo{Proof} By the definition of $B_{s,\vx}$, one has either $
r_{s,\vx} = 2r_0$ (then $B_{s,\vx}\supset B$ and the claim follows from (A2)), or 
$
\|\psi_s\|_{\sssize B(\vx,r)} > \rho  \quad\forall\,r > r_{s,\vx}\,,$ in which case  we are done by the continuity of $\psi_s$. 
\qed\enddemo

For any $\vx\in E$ choose an element $s_\vx$ of $ H(\vx)$ such that   $r_{s_\vx,\vx} \ge  
r_{s,\vx}$ for all $s\in H(\vx)$ (this can be done since $H(\vx)$ is finite). For brevity we will denote $r_{s_\vx,\vx}$ by $r_\vx$ and $B_{s_\vx,\vx}$ by $B_\vx$. Note that in fact $B_\vx = \bigcup_{s\in H(\vx)}B_{s,\vx}$ and $r_\vx = \max_{s\in H(\vx)}r_{s,\vx}$; in particular, 
$
r_\vx \le 2r_0$ $\text{for any }\vx\in E
$.

\proclaim{Lemma 4.4}   For any $\vx\in E$ and $s\in S$, one has 
$\|\psi_s\|_{\sssize B_{\vx}} \ge \rho$.
\endproclaim

\demo{Proof} Assume that $\|\psi_s\|_{\sssize B_{\vx}} < \rho$. 
Then one necessarily has $|\psi_s(\vx)| < \rho $, therefore $s\in H(\vx)$ and  $B_{s,\vx}$ is defined. But $B_{s,\vx}$ is contained in $ B_{\vx}$, so the claim follows from Lemma 4.3.\qed\enddemo

For any $\vx\in E$ we let $S_\vx\df S(\{s_\vx\})$.
The induction procedure of the proof of Theorem 4.1 will be based on the following

\proclaim{Lemma 4.5}   For any $\vx\in E$, one has 
$$
(S_\vx,\ph|_{S_\vx},B_\vx)\in\Cal A(d,k-1,C,\alpha,\rho)\,.
$$ 
\endproclaim

\demo{Proof} The properties (A0) and (A2) for 
$(S_\vx,\ph|_{S_\vx},B_\vx)$ follow from (4.1) and Lemma 4.4 respectively. To prove (A1) and (A3) it suffices to notice that 
$$
\split
\tilde B_\vx\df B(\vx,3^{k-1}r_\vx)&\subset B(\vx_0,3^{k-1}r_\vx + r_0)\\ &\subset B\big(\vx_0,(2\cdot 3^{k-1} + 1)r_0\big)\subset B(\vx_0,3^kr_0) = \tilde B\,.\ \qed
\endsplit
$$
\enddemo 

 The next lemma gives one a way to prove that a point $\vz\in B$ is $(\vre,S,\ph)$-marked provided it is $(\vre,S_\vx,\ph|_{S_\vx})$-marked for some $\vx\in E$.

\proclaim{Lemma 4.6} For  $\vre > 0$ and $\vx\in E$, let $\vz\in B\cap \p(\vre,S_\vx,\ph|_{S_\vx},B_\vx)$ be such that $|\psi_{s_\vx}(\vz)| \ge\vre$.
Then $\vz\in \p(\vre,S,\ph,B)$. Equivalently, 
$$
B\cap \big(B_\vx\ssm\p(\vre,S,\ph,B)\big) \subset B_\vx\ssm\p(\vre,S_\vx,\ph|_{S_\vx},B_\vx) \cup \{\vz\in B_\vx\bigm| |\psi_{s_\vx}(\vz)|  < \vre\}\,.
$$
\endproclaim

\demo{Proof} By the definition of $\p(\vre,S_\vx,\ph|_{S_\vx},B_\vx)$, there exists 
a linearly ordered subset $\Sigma_{\vx,\vz}$ of $S_\vx$ such that $$
\vre\le |\psi_s(\vz)| \le \rho\quad \forall\,s\in\Sigma_{\vx,\vz}\tag 4.2
$$
and
$$
|\psi_s(\vz)| \ge  \rho \quad\forall\,s\in S_\vx(\Sigma_{\vx,\vz})\,.\tag 4.3
$$

Put $\Sigma_{\vz}\df\Sigma_{\vx,\vz}\cup\{s_\vx\}$. Then $S(\Sigma_{\vz}) = S_\vx(\Sigma_{\vx,\vz})$; therefore (M2) immediately follows from (4.3), and, in view of (4.2),  it remains to check (M1) for $s = s_\vx$. The latter is straightforward: $|\psi_{s_\vx}(\vz)|$ is not less than $\vre$ by the assumption and 
is not greater than $\rho $ since $\vz\in B_\vx$. \qed\enddemo

We are now ready to give a 

\demo{Proof of Theorem 4.1} We proceed by induction on $k$.  
First take $\sfb\in \Cal A(d,0,C,\alpha,\rho)$; since the poset $S$ is empty, all points of $B$ are $(\vre,S,\ph)$-marked for  any  $\vre > 0$, which means that  in the case $k=0$ the claim is trivial.  Now take $j\ge 1$ and  suppose that the theorem is proved for $k = j-1$; put  $k = j$ and take $\sfb\in\Cal A(d,k,C,\alpha,\rho)$  and a positive $\vre$.

In view of Lemma 4.2, it suffices to estimate the measure of $E\ssm \p(\vre,S,\ph,B)$. On the other hand, from Lemma 4.6 one deduces that 
 for any $\vx\in E$, the measure of the intersection of $B_\vx\ssm\p(\vre,S,\ph,B)$ with $B$ is not greater than 
$$
\left|B_\vx\ssm\p(\vre,S_\vx,\ph|_{S_\vx},B_\vx)\right|  + \left|\{\vz\in B_\vx\bigm| |\psi_{s_\vx}(\vz)|  < \vre\}\right| \,.
$$
The first summand is not greater than 
$
(k-1)C\big(3^dN_d\big)^{k-1}\left(\dfrac{\vre}\rho \right)^\alpha |B_\vx|
$
 by Lemma 4.5 and the induction assumption. The second one, in view of $\psi_{s_\vx}$ being \cag\ on $\tilde B\supset B_\vx$, is not greater than
$$
C\left(\frac\vre{\|\psi_{s_\vx}\|_{\sssize B_\vx}}\right)^\alpha  |B_\vx|\underset{\text{Lemma 4.4}}\to\le C \left(\frac{\vre} \rho \right)^\alpha  |B_\vx|\,.
$$
Consequently, one gets
$$
\aligned
\left|B\cap \big(B_\vx\ssm\p(\vre,S,\ph,B)\big) \right| &\le C\left((k-1)\big(3^dN_d\big)^{k-1} + 1 \right)\left(\frac{\vre} \rho \right)^\alpha  |B_\vx|\\
&\le kC\big(3^dN_d\big)^{k-1}\left(\frac{\vre} \rho \right)^\alpha  |B_\vx|\,.
\endaligned\tag 4.4
$$

Now consider the covering $\{B_\vx\mid \vx\in E\}$ of $E$ and, using \be's Covering Theorem,  choose a subcovering $\{B_i\}$ of multiplicity  at most $N_{d}$. Then 
$$
\sum_i|B_i| \le N_{d}\left|{\tsize\bigcup_i B_i}\right| \le N_{d}\left|B(\vx_0,3r_0)\right| \le N_{d}3^d|B|\,.\tag 4.5
$$
The sets $B\cap \big(B_i\ssm\p(\vre,S,\ph,B)\big) $ cover $E\ssm \p(\vre,S,\ph,B)$, therefore
$$
\split
\left|\big(E\ssm \p(\vre,S,\ph,B)\right|&\le \sum_i \left|B\cap \big(B_i\ssm\p(\vre,S,\ph,B)\big)\right|\\
\text{(by (4.4))}\quad&\le  kC\big(3^dN_d\big)^{k-1}\left(\frac{\vre} \rho \right)^\alpha \sum_i|B_i|\\ 
\text{(by (4.5))}\quad&\le kC\big(3^dN_d\big)^{k-1}\left(\frac{\vre} \rho \right)^\alpha  N_{d}3^d|B|
\le  kC\big(3^dN_d\big)^k \left(\frac\vre \rho \right)^\alpha |B|\,,
\endsplit
$$
and the theorem is proven. \qed\enddemo

\heading{5. Quantitative non-divergence in the space of lattices}
\endheading 

We now apply Theorem 4.1 to the poset of discrete subgroups of  $\bz^k$, $k\in \bn$. We fix a basis $\ve_1,\dots,\ve_k$ of $\br^k$, and for $I = \{i_1,\dots,i_j\}\subset \{1,\dots,k\}$, $i_1 < i_2 < \dots < i_j$, we let $\ve_{\sssize I} \df \ve_{i_1}\wedge\dots\wedge \ve_{i_j}\subset \bigwedge^j(\br^k)$, with the convention $\ve_\vrn = 1$. We extend the norm $\|\cdot\|$ from   $\br^k$  to the exterior algebra $\bigwedge(\br^k)$ by $\|\sum_{I\subset \{1,\dots,k\}}w_{\sssize I}\ve_{\sssize I}\| = \max_{I\subset \{1,\dots,k\}}|w_{\sssize I}|$.

For a  discrete subgroup $\Gamma$  of $\br^k$, we denote by $\Gamma_\br$ the minimal linear subspace of $\br^k$ containing $\Gamma$. Let  $j =  \text{dim}(\Gamma_\br)$; say that $\vw\in \bigwedge^j(\br^k)$ {\it represents\/} $\Gamma$ if 
$$
\vw =  \cases &1\qquad\qquad\quad\text{\ \ \ if } j = 0\\ &\vv_{1}\wedge\dots\wedge \vv_{j}\quad\text{if }j > 0 \text{ and }\vv_{1},\dots, \vv_{j} \text{ is a basis of }\Gamma\,.\endcases
$$
Clearly the element representing $\Gamma$ is defined up to a sign. Therefore it makes sense to define the norm of $\Gamma$ by 
$
\|\Gamma\| \df \|\vw\| $, where $\vw\text{ represents }\Gamma$.
Note that if $\Gamma$ is a lattice,  the ratio of $
\|\Gamma\|$ and the volume of the quotient space $\br^k/\Gamma$ is uniformly bounded between two positive constants. 

We need the following simple lemma:  

\proclaim{Lemma 5.1} Let $\Gamma$ be a discrete subgroup  of $\br^k$, $\vv\in \br^k \ssm \Gamma_\br$, and let $\Lambda\subset\br^k$ be a discrete subgroup containing both   $\Gamma$ and $\vv$ such that $\Lambda_\br = \Gamma_\br + \br \vv$. Then $\|\Lambda\|\le k\|\Gamma\|\cdot\|\vv\|$, or, equivalently, $\|\vv\| \ge {\|\Lambda\|} / k{\|\Gamma\|}$.
\endproclaim

\demo{Proof} Let $\vw$ represent $\Gamma$; since $\|\Lambda\|$ is not
greater than the norm of $\vw\wedge\vv$, it suffices to show that
$\|\vw\wedge\vv\|\le k\|\vw\|\cdot\|\vv\|$. Write $\vw = \sum_{I\subset \{1,\dots,k\}}w_{\sssize I}\ve_{\sssize I}$ and $\vv = \sum_{i = 1}^kv_i\ve_i$, then 
$$
\split
\|\sum_{I\subset \{1,\dots,k\}}w_{\sssize I}\ve_{\sssize I}\,{\wedge}\, \sum_{i = 1}^kv_i\ve_i\| &\le k\cdot\max\Sb{1\le i \le k}\\{I\subset \{1,\dots,k\}}\endSb|w_{\sssize I}v_i| \\&\le k\cdot\max_{I\subset \{1,\dots,k\}}|w_{\sssize I}|\cdot\max_{1\le i \le k}|v_i| = k\|\vw\|\cdot\|\vv\|\,. \qed
\endsplit
$$
\enddemo

Let $\Lambda$ be a discrete subgroup  of $\br^k$. We say that a subgroup $\Gamma$ of $\Lambda$ is {\it primitive\/} (in $\Lambda$) if $\Gamma = \Gamma_\br\cap \Lambda$, and denote by $\Cal L(\Lambda)$ the set of all nonzero  primitive subgroups of $\Lambda$. The inclusion relation makes $\Cal L(\Lambda)$ a poset, its length being equal to the dimension of $\Lambda_\br$.

 \proclaim{Theorem 5.2} Let $d,k\in\bn$, $C,\alpha > 0$, $0 < \rho  \le 1/k$,  and let a ball $B = B(\vx_0,r_0)\subset \br^d$ and a map $h:\tilde B \to GL_k(\br)$ be given, where $\tilde B$ stands for $B\big(\vx_0,3^kr_0\big)$.  
 For any $\Gamma \in \Cal L(\bz^k)$, denote by $\psi_{\sssize \Gamma}$ the function $\psi_{\sssize \Gamma}(\vx) \df \|h(\vx)\Gamma\|$, $\vx\in \tilde B$. Assume that for any $\Gamma \in \Cal L(\bz^k)$,
\roster
\item"(i)" $\psi_{\sssize \Gamma}$ is \cag\ on $\tilde B$;
\item"(ii)" $\|\psi_{\sssize \Gamma}\|_{\sssize B} \ge \rho $.
\endroster
Then 
 for any  positive $ \vre \le \rho$ one has
$$
\left|\{\vx\in B\mid \delta\big(h(\vx)\bz^k\big) < \vre\}\right| \le kC \big(3^dN_d\big)^k \left(\frac\vre \rho \right)^\alpha  |B|\,.\tag 5.1
$$
\endproclaim

\demo{Proof} We let $S = \Cal L(\bz^k)$ and denote by $\ph$ the map $\Gamma\to \psi_{\sssize \Gamma}$. It is easy to verify that  $\sfb\in\ca$.
Indeed, the property (A0) is clear, (A1) is given by (i), (A2) by (ii), and (A3) follows from the discreteness of $\bigwedge(\bz^k)$ in $\bigwedge(\br^k)$. 

In view of Theorem 4.1, it remains to prove that a point $\vx\in B$ with  $\delta\big(h(\vx)\bz^k\big) < \vre$ can not be  $(\vre,S,\ph)$-marked. In other words,
$$
\p(\vre,S,\ph,B) \subset  \{\vx\in B\mid \delta\big(h(\vx)\bz^k\big) \ge \vre\}\,.\tag 5.2
$$

Take an $(\vre,S,\ph)$-marked point $\vx\in B$, and let $\{0\} = \Gamma_0 \subsetneq \Gamma_1 \subsetneq\dots \subsetneq\Gamma_l = \bz^k$ be all the elements of $\Sigma_\vx \cup \big\{\{0\},\bz^k\big\}$. Take any $\vv\in\bz^k\nz$. Then there exists $i$, $1 \le i \le l$, such that $\vv\in \Gamma_i\ssm \Gamma_{i-1}$. Denote $(\Gamma_{i-1} + \br\vv)\cap \bz^k$ by $\Lambda$. Clearly $\Lambda$ is a primitive subgroup of $\bz^k$ contained in $\Gamma_i$, therefore $\Lambda \in \Sigma_\vx \cup S(\Sigma_\vx)$. Now one can use properties (M1) and (M2) to deduce that
$$
|\psi_{\sssize \Lambda}(\vx)| = \|h(\vx)\Lambda\| \ge \min(\vre, \rho) = \vre\,,
$$
and then apply Lemma 5.1 to conclude that 
$$
\|h(\vx)\vv\| \ge \frac{\|h(\vx)\Lambda\|}{k\|h(\vx)\Gamma_{i-1}\|} \ge \vre/k\rho  \ge \vre\,.
$$
This shows (5.2) and completes the proof of the theorem. \qed\enddemo

As was mentioned in the introduction, our method of proof is based, with some technical changes, on the argument from \cite{Mar} and its modification in \cite{D2}. As an illustration, let us show how one can
use the above theorem to get a quantitative strengthening of Theorem 2.1 from \cite{D2}.

 \proclaim{Theorem 5.3} For any lattice  $\Lambda$ in $\br^k$ there exists a constant $\rho = \rho(\Lambda) > 0$ such that for any one-parameter unipotent subgroup $\{u_x\}_{x\in\br}$   of $SL_k(\br)$, for any $T > 0$ and any $\vre \le \rho$, one has
$$
\left|\{0 < x < T\mid \delta\big(u_x\Lambda\big) < \vre\}\right| \le
2k^36^k (k^2 + 1)^{1/k^2} \left(\frac\vre \rho \right)^{1/k^2}T\,.\tag
5.3 
$$
\endproclaim

\demo{Proof} Write $\Lambda$ in the form  $g\bz^k$ with $g\in GL_k(\br)$, and denote by $h$ the function $h(x) = u_xg$. For any $\Gamma \in \Cal L(\bz^k)$ with basis $\vv_{1},\dots, \vv_{j}$, the coordinates of \linebreak $h(x)(\vv_{1}\wedge\dots\wedge \vv_{j})$ will be polynomials in $x$ of degree not exceeding $k^2$, hence  \linebreak $\big(2k^2 (k^2 + 1)^{1/k^2},{1/k^2}\big)$-good on $\br$ by Proposition 3.2. In view of  Lemma 3.1(c) and the definition of the norm of a lattice, the functions
$\psi_{\sssize \Gamma}(x) \df \|h(x)\Gamma\|$ will also be  $\big(2k^2 (k^2 + 1)^{1/k^2},{1/k^2}\big)$-good on $\br$. 

Now let $\rho \df \min\big(1/k,\inf_{\Gamma \in \Cal L(\bz^k)}\|g\Gamma\|\big)$, positive by the discreteness of $\Lambda$ in $\br^k$. Then $\psi_{\sssize \Gamma}(0) \ge \rho$ for any $\Gamma \in \Cal L(\bz^k)$, therefore $\|\psi_{\sssize \Gamma}\|_{(0,T)}\ge \rho$ by the continuity of  $h$.

We see now that with the the substitutions $B = (0,T)$, $C = 2k^2
(k^2 + 1)^{1/k^2}$, $\alpha = {1/k^2}$ and $d = 1$ (note that it is an
elementary fact that $N_1 = 2$), assumptions (i) and (ii) of Theorem
5.2 are satisfied, and one immediately gets (5.3) from
(5.1). \qed\enddemo 

We now derive another corollary from Theorem 5.2 which will immediately imply Proposition 2.3. In what follows, we put $k = n + 1$ and for $\vy\in\br^n$ define $u_\vy \df \left(\matrix
1 & \vy^{\sssize T}  \\
0 & I_n
\endmatrix \right)$, so that the lattice $\Lambda_\vy$ (see (1.4)) is given by  $u_\vy \bz^{n+1}$. We also let $\{\ve_0,\ve_1,\dots,\ve_n\}$ be the standard basis of $\br^{n+1}$.

 \proclaim{Theorem 5.4}  For $n\in \bn$,  put $k = n + 1$ and let $d$,
$C$, $\alpha$, $\rho$ and  
$B = B(\vx_0,r_0)$ be as in Theorem 5.2.  Also let $\vf = (f_1,\dots,f_n)$ be a continuous map from $\tilde B \df B\big(\vx_0,3^{n+1}r\big)$ to $\br^n$ such that
\roster
\item"(i)" for any $\vc = (c_0,c_1,\dots,c_n) \in \br^{n+1}$, $c_0 + \sum_{i = 1}^nc_if_i$ is \cag\ on $\tilde B$; 
\item"(ii)" for any $\vc\in \br^{n+1}$ with $\|\vc\| \ge 1$, $\|c_0 + \sum_{i = 1}^nc_if_i\|_{\sssize B} \ge \rho $.
\endroster
Fix nonnegative numbers $t_1,\dots,t_n$,  
and for  $\vx\in \tilde B$ let
$$
h(\vx) \df  g_\vt u_{\vf(\vx)}\,,\tag 5.4
$$
where $g_\vt $ is given by {\rm (1.5)}. 
Then 
 for any  positive $ \vre \le \rho$ one has
$$
\left|\{\vx\in B\mid \delta\big(h(\vx)\bz^{n+1}\big) < \vre\}\right| \le ({n+1})C\big(3^dN_d\big)^{n+1} \left(\frac\vre \rho \right)^\alpha  |B|\,.\tag 5.5
$$
\endproclaim


\demo{Proof} 	In view of Theorem 5.2, it suffices to show that (i) and (ii) above will remain true with $c_0 + \sum_{i = 1}^nc_if_i$ replaced with $\psi_{\sssize \Gamma} \df \|h(\cdot)\Gamma\|$, where $\Gamma$ is any element of  $\Cal L(\bz^{n+1})$. Take   $\vw = \sum_{I\subset \{0,\dots,n\}}w_{\sssize I}\ve_{\sssize I}$ representing $\Gamma$. To see how the coordinates of $\vw$ change under the action of $h(\vx)$ of the form (5.4), first note that the action of $u_{\vf(\vx)}$ leaves $\ve_0$ invariant and sends $\ve_{i}$ to $\ve_{i} + f_i(\vx) \ve_0$, $i = 1,\dots,n$. Therefore 
$$
u_{\vf(\vx)}\ve_{\sssize I} =  \cases &\ve_{\sssize I} \text{ if }0\in I\\
&\ve_{\sssize I} + \sum_{i\in I} \pm f_i(\vx)
\ve_{\sssize I \cup \{0\}\ssm\{i\}}\text{ otherwise}\,,\endcases
$$
which shows that
$$
u_{\vf(\vx)}\vw =  \sum_{0\notin I}w_{\sssize I}\ve_{\sssize I}
+ \sum_{0\in I}\left( w_{\sssize I} + \sum_{i\notin I} \pm w_{\sssize I\cup\{i\}\ssm\{0\}}f_i(\vx) \right) \ve_{\sssize I}\,.
$$
Applying $g_\vt$ to both sides of the above formula, one gets $
h(\vx)\vw = \sum_{I\subset \{0,\dots,n\}}h_{\sssize I}(\vx)\ve_{\sssize I}$, where
$$
h_{\sssize I}(\vx) = \cases & e^{-\sum_{i\in I}t_i} w_{\sssize I} \text{ if }0\notin I\\
&e^{\sum_{i\notin I}t_i}\left(w_{\sssize I} + \sum_{i\notin I} \pm w_{\sssize I\cup\{i\}\ssm\{0\}}f_i(\vx) \right) \text{ otherwise}\,.\endcases
$$

It follows that all the coordinates $h_{\sssize I}(\vx)$ of $h(\vx)\vw$ are of the form $c_0 + \sum_{i = 1}^nc_if_i(\vx)$ for some $\vc \in \br^{n+1}$. Thus, by virtue of assumption (i) of the Theorem, they are  \cag\ on $\tilde B$. Then Lemma 3.1(c) applies and one can conclude that $\psi_{\sssize \Gamma} = \sup_{I\subset \{0,\dots,n\}}|h_{\sssize I}|$ is \cag\ on $\tilde B$. Furthermore, since the coordinates $w_{\sssize I}$ of $\vw$ are integers and at least one of them is nonzero, one can conclude that for some $I\ni 0$, the function   $h_{\sssize I}(\vx)$ is of the form $c_0 + \sum_{i = 1}^nc_if_i$ with  $\|\vc\| \ge 1$. Hence $\|h_{\sssize I}\|_{\sssize B} \ge \rho$ for this $I$. Therefore $\|\psi_{\sssize \Gamma} \|_{\sssize B} \ge \rho$, which is all one needs to apply Theorem 5.2. 
\qed\enddemo

It is now easy to write down the

\demo{Proof of Proposition 2.3} Take $U\subset \br^d$, $\vf:U\to\br^n$, $\vx_0\in U$ and $l$ as in the statement of the proposition. Using Proposition 3.4, 
find a neighborhood $V\subset U$ of $\vx_0$ and $C> 0$  such that any linear combination of $1,f_1,\dots,f_n$ is $(C,1/dl)$-good on $V$. Choose a ball $\tilde B = B(\vx_0,\tilde r)$ contained in $V$ and let $B = B\big(\vx_0,3^{-(n+1)}\tilde r\big)$. Then  
condition (i) of Theorem 5.4 is satisfied with $\alpha = 1/dl$, while
the existence of positive $\rho$ satisfying (ii) follows from the linear independence of $1,f_1,\dots,f_n$ on $B$ over $\br$. The validity of (2.6), with $D = ({n+1})C\big(3^dN_d\big)^{n+1}$, is now an immediate consequence of (5.5).
\qed\enddemo

\heading{6. Concluding remarks and generalizations}
\endheading

\subhead{6.1}\endsubhead 
Observe that a notion of \cag\ functions can be defined for any metric space $X$ with a Borel measure $\mu$ in place of $(\br^d,|\cdot|)$. Moreover, the argument of \S 4 works for arbitrary $(X,\mu)$ provided two additional conditions are satisfied:

\roster 
\item"{$\bullet$}" Besicovitch's Covering Theorem holds, the constant $N_d$ being replaced by some positive number $N_{\sssize X}$;

\item"{$\bullet$}" (cf.~(4.5)) $\dsize\sup_{x\in X,\,r > 0} \ \frac{\mu\big(B(x,3r)\big)}{\mu\big(B(x,r)\big)} < \infty$.

\endroster

These conditions, in particular, are satisfied for $X = \bq_p$ and $\mu$ a Haar measure on $\bq_p$. This way one can attempt to apply the methods of the present paper to describe extremal and strongly extremal manifolds in $p$-adic spaces, generalizing the results from \cite{Sp2, Chapter 2}. 
This work is currently in progress.

\subhead{6.2}\endsubhead One may ask whether it is possible
to develop a similar proof starting from the inequality (1.2M)
instead of (1.1M). The answer is yes, although the proof turns
out to be slightly more complicated. The recipe is simple: one should consider
lattices
$
\left(\matrix
I_n & \vy  \\
0 & 1
\endmatrix \right)\bz^{n+1}
$ instead of (1.4) and act on them by elements of the form \linebreak 
$\text{\rm diag}(e^{t_1},\dots,e^{t_n}, e^{-t})$ instead of (1.5). 

More generally, one can unify these two approaches by saying that
a matrix  $Y\in \text{Mat}_{m,n}(\br)$, interpreted as a system   of $m$ linear forms in $n$ variables, is  VWA (resp.~VWMA)  if for some $\vre > 0$ there are infinitely many $\vq\in \bz^n$ such that 
$$
\|Y\vq + \vp\|^m\cdot\|\vq\|^n \le \|\vq\|^{-n\vre}\,,\tag 6.1
$$
or, respectively,
$$
\Pi(Y\vq + \vp)\cdot\Pi_{\sssize +}(\vq) \le \Pi_{\sssize +}(\vq)^{-\vre}\,,\tag 6.1M
$$
for some $\vp\in\bz^m$. Then one faces a problem of describing (strongly) extremal submanifolds of $\text{Mat}_{m,n}(\br)$ \footnote{See \cite{Ko1, Ko2} where similar objects are called ``systems of (strongly) jointly extremal manifolds".}. It turns out that one can apply Theorem 5.2 to get a general result of which Theorem A is a special case. This is going to be a topic of a forthcoming paper.

\subhead{6.3}\endsubhead It is instructive to compare Theorem A with results recently obtained by M.~Dodson, B.~Rynne and J.~Vickers. Following \cite{DRV2}, say that the manifold $M$ satisfies condition $K1$ at $\vy\in M$ if for any $\vv\in T_\vy M^\perp$, at least two of the principal curvatures of $M$ at $\vy $ with respect to $\vv$ (see \cite{DRV1} for further details) are nonzero and have the same sign. It is proved in \cite{DRV1, DRV2} that a $C^3$ manifold $M$ is extremal provided $K1$ holds for almost all $\vy\in M$. Now say that $M$ satisfies condition $K0$ at $\vy\in M$ if for any $\vv\in T_\vy M^\perp$, at least one  principal curvature of $M$ at $\vy $ with respect to $\vv$ is nonzero. Then from Theorem A it immediately follows that
a $C^2$ submanifold $M$ of $\br^n$ is strongly extremal whenever
 condition $K0$ holds for almost every its point.
Indeed, using the analytic formulation of curvature conditions as in \cite{DRV1}, it is easy to see that $M$ satisfies condition $K0$ at $\vy$ if and only if $\vy$ is a $2$-nondegenerate point of $M$.

On the other hand, Dodson, Rynne and Vickers were able to use stronger curvature conditions to derive several Khintchine-type theorems and asymptotic formulae, see \cite{DRV2, DRV3}. It would be interesting to know whether it is possible to obtain similar results assuming only the mild condition $K0$. In particular, it is proven in \cite{DRV2} that if $K1$ holds for almost all $\vy\in M$ and $\ph:\bn\to\br_{\sssize +}$ is a decreasing function such that
$$
\sum_{q = 1}^\infty \frac{\ph(q)}q  < \infty\,,\tag 6.2
$$
 then for a.e.~$\vy\in M$ the inequality
$$
|\vq\cdot\vy + p|\cdot\|\vq\|^n \le \ph\big(\|\vq\|^{n}\big)\tag 6.3
$$
has at most finitely many solutions $\vq\in\bz^n$, $p\in \bz$. On the other hand, one can modify Lemma 2.1  in the spirit of \cite{K1, Theorem 8.3}, and then apply Proposition 2.3 with $\vt = (t/n,\dots,t/n)$  to obtain the following 

\proclaim{Theorem B} Let $f_1,\dots,f_n\in C^l(U)$, $U$ an open subset of $\br^d$,  be such that a.e.~point of $M = \{\vf(\vx)\mid \vx\in U\}$ is {$l$-nondegenerate}, and let $\ph:\bn\to\br_{\sssize +}$ be a decreasing function such that
$$
\sum_{q = 1}^\infty \frac{\ph(q)^{1/dl(n+1)}}q  < \infty\,.\tag 6.4
$$
 Then for a.e.~$\vy\in M$ the inequality {\rm (6.3)}
has at most finitely many solutions.
\endproclaim

The details will appear elsewhere; whether it is possible to replace (6.4) by (6.2) in the above theorem is an 
 open question. We remark that for $M$ of the form (1.3) and $\ph$ satisfying (6.2), the conclusion of Theorem B was conjectured by A.~Baker in 1966 and proved by V.~Bernik \cite{Be} in 1984.


\heading{Acknowledgements}
\endheading 

A substantial part of this work was done during the authors' stay at the University of Bielefeld in June--July 1996. This stay was supported by SFB-343 and Humboldt Foundation.
Thanks are also due to V.~Bernik, D.~Dickinson, M.~Dodson,
T.~Khovanova, N.~Shah, G.~Tomanov and the referee for useful comments.

\enddocument